\documentclass[11pt]{amsart}

\usepackage{amsmath,amssymb}
\usepackage{amsfonts}
\usepackage{amsthm}
\usepackage{latexsym}
\usepackage{graphicx}
\usepackage{epsfig}

\usepackage{color}

\def\ppt{\frac{\partial}{\partial t}}

\numberwithin{equation}{section}
\newtheorem{theorem}{Theorem}
\newtheorem{definition}{Definition}

\newtheorem{lemma}{Lemma}

\newtheorem{remark}{Remark}
\numberwithin{theorem}{section}
\numberwithin{definition}{section}
\numberwithin{lemma}{section}
\numberwithin{proposition}{section}
\numberwithin{remark}{section}
\def\R{\bf R}

\def\al{\aligned}
\def\eal{\endaligned}
\def\M{{\bf M}}
\def\be{\begin{equation}}
\def\ee{\end{equation}}
\def\lab{\label}

\begin{document}

\tracingpages 1
\title[the conjugate heat equation and Ancient solutions of the Ricci flow]
{\bf the
conjugate heat equation and Ancient solutions of the Ricci flow}

%    Information for first author
\author{Xiaodong Cao $^*$}
\thanks{$^*$Research partially supported by NSF grant DMS 0904432 and by the Jeffrey Sean Lehman Fund from Cornell University}

%    Address of record for the research reported here
\address{Department of Mathematics,
  Cornell University, Ithaca, NY 14853-4201}
\email{cao@math.cornell.edu}
%    Current address
%\curraddr{}
%    \thanks will become a 1st page footnote.
%\thanks{}

%    Information for second author
\author{Qi S. Zhang}
%    Address of record for the research reported here
\address{Department of Mathematics, University of California, Riverside, CA
92521} \email{qizhang@math.ucr.edu}
%    Current address
%\curraddr{}
%    \thanks will become a 1st page footnote.
%\thanks{}

\renewcommand{\subjclassname}{%
  \textup{2000} Mathematics Subject Classification}

\date{May 28, 2010}

\begin{abstract}
We prove Gaussian type bounds for the fundamental solution of the conjugate heat equation evolving under the Ricci flow. As a consequence, for dimension $4$ and higher, we show that the backward limit of type I $\kappa$-solutions of the Ricci flow must be a non-flat gradient shrinking Ricci soliton.  This extends Perelman's  previous result on backward limits of $\kappa$-solutions in dimension $3$, in which case that the curvature operator is nonnegative (follows from Hamilton-Ivey curvature pinching estimate). The Gaussian bounds that we obtain on the fundamental solution of the conjugate heat equation under evolving metric might be of independent interest.

\end{abstract}
\maketitle \tableofcontents
\section{Introduction}

Let $(\M, g(t))$ be a complete solution with bounded curvature to the Ricci flow
\be \label{rf}
\ppt g(t) = - 2 Ric_{g(t)}.
\ee
 It is called an ancient solution if it is defined for all $t\in (-\infty, T_0)$, for some $T_0>0$. Ancient solutions typically arise as singularity models of the Ricci flow. For example,  Type I singularity model is a dilation limit of a Type I maximal solution to the Ricci flow. As a consequence, classifications of ancient solution are very important subjects in the study of the Ricci flow.
It is well-known that all ancient solutions have nonnegative
scalar curvature (for example, see \cite{chenbl09}).\\

In \cite{perelman1}, G. Perelman showed that the rescaling limits at
singularities of the Ricci flow are $\kappa$-solutions, which is
defined as the following,

\begin{definition}
\lab{defgujie} A complete, non-flat ancient solution $({\bf M},
g(t))$, $t\in (-\infty, T_0)$,  $T_0 > 0$, to the Ricci flow is a
$\kappa$-solution if it  is $\kappa$-non-collapsed on all scales
for some positive constant $\kappa$, i.e., $\forall (x_0,t_0) \in
{\bf M}\times (-\infty, T_0)$, $\forall r>0$,
 let $P(x_0, t_0, r, -r^2)$ be the  parabolic ball
 \[
 \{(x, t) \ | \
 d(x, x_0, t) < r, \quad t_0-r^2 < t < t_0 \ \},
 \]
 and if  $|Rm| \le r^{-2}$ on $P(x_0, t_0, r, -r^2)$, then we have
$|B(x_0, r, t_0)|_{t_0} \ge  \kappa r^n$, here $|B(x_0, r, t_0)|_{t_0}$ stands for the volume of the ball centered at $x_0$ with radius $r$ at time $t_0$, which is measured using the metric $g(t_0)$.
%${\bf M}$ is $\kappa$-non-collapsed at $(x_0, t_0)$ at scale $r$
\end{definition}
\medskip

In dimensional $3$, Perelman proves that all $\kappa$-solutions
have a rescaled backward limit in time, which is a non-flat
gradient shrinking Ricci soliton. By classifications of
gradient shrinking Ricci solitons of T. Ivey \cite{Isoliton} (for
compact case, also see \cite{perelman1}) and Perelman
\cite{perelman2} (for complete non-compact case), Perelman was
able to obtain some qualitative result about $\kappa$-solutions.
Moreover, this leads to properties of canonical neighborhoods for
Ricci flows, which is important in the study of Ricci flows in
dimension $3$. Therefore, we believe that, the understanding of
backward limits of $\kappa$-solutions will also play an important
role in the study of the singularity of Ricci flow in high
dimensions. One hopes that all such backward limits still remain
as gradient shrinking Ricci
solitons.\\

In this paper, we generalize Perelman's classification result on backward limit solution, namely, we prove that, in all dimensions, for a non-flat, Type I $\kappa$-solution of the Ricci flow, there exists a sequence of backward dilation solutions, which converge to a non-flat gradient shrinking Ricci soliton. We first recall the definition of a Type I solution,

\begin{definition}  A $\kappa$-solution on $(-\infty, T_0)$ is called Type I if there exists
a positive constant $D_0$ such that $$|Rm(x, t)|
\le \frac{D_0}{T_0-t}.$$
\end{definition}
\noindent Notice that here we only require that the solution to be Type I, but not assuming that the curvature operator of such  $\kappa$-solutions to be nonnegative (which is the case of previous results for dimension $4$ and higher). In \cite{perelman1}, Perelman introduced reduced distance and reduced volume to classify the backward limits in dimension $3$. Nonnegativity of the curvature operator (or sectional curvature) plays an essential role in that proof. Our approach here is different, the main tools that we use in this paper are the W-entropy (also introduced by Perelman) and bounds on the fundamental solution for the conjugate heat equation associated with the Ricci flow (\ref{rf}),
\be
\lab{conjheat}
%\begin{cases}
\ppt u=-\Delta u + R u, % \\
%\partial_t g(t) = - 2 Ric_{g(t)}.
%\end{cases}
\ee here $\Delta$ is the Laplace-Beltrami operator,   $R$ is the scalar curvature, both are  with respect to the metric $g(t)$.

The original idea that to study the Ricci flow coupled with the Harmonic Map flow arose from R. Hamilton's paper \cite{Hsurvey}.  Perelman successfully proved a non-collapsing result using the conjugate heat equation associated with the Ricci flow, the authors (\cite{c08} and \cite{kz08}) proved a differential Harnack inequality for positive solutions of (\ref{conjheat}). The system of the heat equation with the Ricci flow was studied by many authors, various estimates and applications can be found in \cite{guenther02, ni04, z06, act08, ch2009, liu09, bcp10}. It turned out that the technical results and additional information one obtained from these associated systems were well worth the extra difficulty in analysis caused by the extra equation(s).

Estimates of the fundamental solution of the heat equation on manifolds has been a traditionally active research area with many applications. In the fixed metric case, we refer the reader to \cite{ly86, davies90} and \cite{grigoryan97} for more information on this vast field. For the conjugate heat equation, existence
of the fundamental solution was first proved by C. Guenther in \cite{guenther02}.  Various interesting bounds were obtained by various authors in \cite{guenther02, ni04, z06, chowetc1} and \cite[Section 9]{perelman1}. We comment here that the fundamental solution of the heat equation can also be viewed as the fundamental solution of the conjugate heat equation in a very natural way, and we will use this simple but important observation from time to time throughout this paper.\\

The rest of this paper is organized as follows. In Sections 2, we first prove an on diagonal upper bound for the fundamental solution. One of the main ingredients that used in the proof is the uniform Sobolev inequality under the Ricci flow (for example, see \cite{z07}). As explained above, what we prove is an upper bound for the fundamental solution of the forward heat equation coupled with the Ricci flow. Then the lower bound of the fundamental solution essentially follows from a differential Harnack inequality proved in \cite{z06} and \cite{ch2009}.

In Section 3, we restrict ourselves to the case of nonnegative Ricci curvature. We establish certain  point-wise bound for the fundamental solution of the conjugate heat equation coupled with the Ricci flow. The result itself can be viewed as a Ricci flow version of \cite{ly86}, where P. Li and S.-T. Yau proved sharp bounds for the heat kernel on manifolds with nonnegative Ricci curvature. However, our technique here is much different due to the fact that the metric is evolving by the Ricci flow. Beside using the classical method of A. Grigor'yan \cite{grigoryan97}, the uniform Sobolev inequality under the Ricci flow (\cite{z07}) and the elliptic type Harnack inequality (\cite{z06, ch2009}) again plays an key role in the proof. The result in this section may has its own interest in analysis and probability, simply as an analogous result for the fixed metric case.

In Section 4, we prove that, the backward limits of type I $\kappa$-solutions is a non-flat gradient shrinking Ricci soliton, using the bounds that we obtained in previous sections. One application of this result is the characterization of Type I singularity model of the Ricci flow. Since our limit might be as well complete non-compact, we first extend out estimate in previous sections to the complete non-compact case. Here we assume that the solution is of Type I, while previous results in high dimensions all need to assume nonnegative curvature operators. Since the Type I condition plays a crucial role in proving the Gaussian bounds and rescaling process, it remains interesting to see if similar result holds for type II $\kappa$-solutions.

Then in Section 5, for the readers' convenience, we state the uniform Sobolev inequality under the Ricci flow, which is used in our proofs.
\\

We now fix our notations here. ${\bf M}$ denotes a complete Riemannian manifold, either non-compact or compact without boundary;  $g_{ij}$ and $R_{ij}$ denote the metric and Ricci curvature in local coordinate systems; $\nabla$ and $\Delta$ denote the corresponding gradient and Laplace-Beltrami operators of the evolving metric $g(t)$;  $d(x, y, t)$ denotes the distance function with respect to $g(t)$; $d\mu_{g(t)}(x)$ denotes the volume element of $g(t)$ at $x$. In the case that ${\bf M}$ is complete non-compact, by fundamental solution we mean the minimal positive fundamental solution to the (conjugate) heat equation. Throughout this paper, we use $c$, $C$, $\alpha$, $\beta$ and $\eta$ (all with or without index) to denote generic positive constant that may change from line to line.\\

We note that there are some recent work on Ricci flow singularities independently by N. Q. Le and N. Sesum \cite{ls2010}, J. Enders, R. M\"{u}ller and P. M. Topping \cite{emt2010}, both of their works study the forward limit of singularities. For Type I solutions, one can take a backward limit of these forward limits, which are $\kappa$-solutions. Notice that limits of singularity models are still singularity models, so our result also implies that Type I singularity model is a non-flat gradient shrinking Ricci soliton.

%%%%%%%%%%%%%%%%%%%%%%%%%%%%%%%%%%%%%%%%%%%%%%%%%%%%%%%%%%%%%%%%%%%%%%%%%%%%%
\section{On diagonal bounds  with no curvature assumptions }
\label{sec:2}
%%%%%%%%%%%%%%%%%%%%%%%%%%%%%%%%%%%%%%%%%%%%%%%%%%%%%%%%%%%%%%%%%%%%%%%%%%%%%

In this section we establish an upper bound for $G=G(x, t; y, s)$
when $x, y$ are close in certain sense. This bound is sometimes
referred to as the on diagonal bound. In the paper \cite{z06},
this kind of bound depending on the Sobolev constants of the
manifolds $({\bf M}, g(t))$ was proven.  In view of the recent
uniform Sobolev inequality under Ricci flow \cite{z07} (also see the
appendix of the paper), we can prove a better upper bound now. In addition we
also prove a corresponding lower bound. These bounds are global in time.\\

In the following, we use $d\mu_{g(t)}$ to denote the volume form
of $g(t)$ at time $t$, and we will omit the subindex $g(t)$ when
there is no confusion. We use $\sup R^-(x, l)$ to denote the
larger number between $-\inf_{\bf M} R(x,l)$ and $0$, i.e., $\sup
R^-(x, l)=\max \{-\inf_{\bf M} R(x,l) ,0\}$. We also use
the short notation $R=R(x,t)$ to denote the scalar curvature at
point $x$ and time $t$ when there is no confusion.  We define $\lambda_0$
be the lowest eigenvalue of  the operator $-4 \Delta +R$ on  $({\bf M}, g(0))$, i.e.,
\begin{equation}
\label{lambda0}
\lambda_0 = \inf_{ \Vert v \Vert_2 = 1} \int_{\M} ( 4 |\nabla v |^2
+ R v^2) d\mu_{g(0)}.
\end{equation}

%\medskip
We first state our main theorem of this section.

\begin{theorem}\label{thm:2.1}
%th2.1
 Let $G=G(z, l; x, t)$, $l<t$,  be the fundamental solution of (\ref{conjheat}).
 Assume that at time $l$, there exist positive constants $A_0, B_0$  such that, for all $v
\in W^{1, 2}({\bf M}, g(l))$, the following Sobolev imbedding theorem holds,
\[
\bigg( \int v^{2n/(n-2)} d\mu_{g(l)} \bigg)^{(n-2)/n} \le A_0 \int
 |\nabla v |^2  d\mu_{g(l)} + B_0 \int v^2
d\mu_{g(l)}.
\]  Then there exist positive constants $c_1$ and $c_2$ such
that
\be
\label{odul}  \frac{c_1 B^{-1}((t-l)/2)}{ (t-l)^{n/2}}
e^{-2 c_2 \frac{d(z, x, t)^2}{t-l}}  e^{- \frac{1}{ \sqrt{t-l}}
\int^{t}_{l} \sqrt{t-s}
 R(x, s) ds} \le G(z, l; x, t) \le \frac{B(t-l)}{(t-l)^{n/2}},
 \ee where $B(t-l)=\exp[\alpha +  (t-l) \beta +
 (t-l) \sup R^-(x, l)]$,  $\alpha=\alpha(A_0, B_0, \lambda_0, n)$ and
 $\beta=\beta(A_0, B_0, \lambda_0, n)$ are positive constants.

 Moreover, in the special case that $R(x, 0) >0$, the function $B(t-l)$ becomes
a constant independent of $t$ and $l$, and (\ref{odul}) becomes
\be \label{odulp}  \frac{c_1B^{-1}}{ (t-l)^{n/2}} e^{-2 c_2
\frac{d(z, x, t)^2}{t-l}}  e^{- \frac{1}{ \sqrt{t-l}} \int^{t}_{l}
\sqrt{t-s}
 R(x, s) ds} \le G(z, l; x, t) \le \frac{B}{(t-l)^{n/2}},
 \ee where $B=\exp[\alpha]$.
\end{theorem}

\begin{remark} (\ref{odulp}) is also true if the scalar curvature is
nonnegative and the solution to the Ricci flow is non-flat, which
is the case of non-flat ancient solutions of the Ricci flow.\end{remark}

\proof Let $G=G(z,l; x, t)$, $l<t$,  be the fundamental solution
of (\ref{conjheat}).  Then as function of $(x, t)$, G is the
fundamental solution of the forward heat equation associated with
the Ricci flow, i.e.,

\begin{equation}\label{forheat}
\begin{cases}\begin{array}{lll}\frac{\partial}{\partial
t}g(t)&=&-2Ric,\\
\ppt u&=&\Delta u.\end{array}\end{cases}\end{equation} The idea
is to study the forward heat equation (\ref{forheat}) first.
Without loss of generality, we may assume that $l=0$.\\

Let $u=u(x,t)$ be a positive solution to (\ref{forheat}). Given $T
>0$ and $t \in (0, T)$, define $p(t)=T/(T-t)$, so $p(0)=1$ and
$p(T)=\infty$. By direct computation
\[
\al
\partial_t \Vert u \Vert_{p(t)} &= \partial_t \bigg[ \bigg(\int_{\bf M}
u^{p(t)}(x, t) d\mu_{g(t)} \bigg)^{1/p(t)} \bigg]\\
&=-\frac{p'(t)}{p^2(t)} \Vert u \Vert_{p(t)} \ln \int_{\bf M}
u^{p(t)}(x, t) d\mu_{g(t)}
    +\frac{1}{p(t)} \bigg( \int_{\bf M} u^{p(t)}(x, t) d\mu_{g(t)} \bigg)^{(1/p(t))-1}
          \\
          &\qquad \times \bigg[ \int_{\bf M} u^{p(t)} (\ln u) p'(t) d\mu_{g(t)} +
    p(t) \int_{\bf M} u^{p(t)-1} ( \Delta u - R u ) d\mu_{g(t)} \bigg].
\eal
\] Using integration by parts on the
term containing $\Delta u$ and multiplying both sides by $p^2(t)
\Vert u \Vert^{p(t)}_{p(t)}$, we arrive at
\[
\al &p^2(t)  \Vert u \Vert^{p(t)}_{p(t)}
\partial_t \Vert u \Vert_{p(t)}\\
=& -p'(t) \Vert u \Vert^{p(t)+1}_{p(t)} \ln \int_{\bf M}
u^{p(t)}(x, t) d\mu_{g(t)}
  + p(t) \Vert u \Vert_{p(t)} p'(t) \int_{\bf M}
      u^{p(t)} \ln u (x, t) d\mu_{g(t)}\\
      &-p^2(t)(p(t)-1) \Vert u \Vert_{p(t)} \int_{\bf M}
            u^{p(t)-2} |\nabla u|^2(x, t) d\mu_{g(t)}
              - p^2(t) \Vert u \Vert_{p(t)} \int_{\bf M}  R(x, t)
            u^{p(t)}(x, t) d\mu_{g(t)}.
\eal
\]Dividing both sides by $ \Vert u \Vert_{p(t)}$, we obtain
\[
\al & p^2(t)  \Vert u \Vert^{p(t)}_{p(t)}
\partial_t \ln \Vert u \Vert_{p(t)}\\
=& -p'(t) \Vert u \Vert^{p(t)}_{p(t)} \ln \int_{\bf M} u^{p(t)}
d\mu_{g(t)}
  + p(t)  p'(t) \int_{\bf M}
      u^{p(t)} \ln u  d\mu_{g(t)}\\
      &-4[p(t)-1]  \int_{\bf M}
            |\nabla (u^{p(t)/2}) |^2 d\mu_{g(t)}
              - p^2(t)  \int_{\bf M}  R
            (u^{p(t)/2})^2 d\mu_{g(t)}.
\eal
\]
Define $v(x,t) = \frac{u^{p(t)/2}}{\Vert u^{p(t)/2} \Vert_2}$, we
have $\Vert v \Vert_2 = 1$ and $$v^2\ln v^2=p(t) v^2\ln u-2v^2 \ln
\Vert u^{p(t)/2} \Vert_2,$$ merging the first two terms on the
righthand side of the above equality and dividing both sides by
$\Vert u \Vert^{p(t)}_{p(t)}$, we arrive at
\[
\al & p^2(t)
\partial_t \ln \Vert u \Vert_{p(t)}\\
=& p'(t)   \int_{\bf M} v^2 \ln v^2  d\mu_{g(t)}
         -4(p(t)-1)  \int_{\bf M}
            |\nabla v |^2 d\mu_{g(t)}
              - p^2(t)  \int_{\bf M}  R v^2 d\mu_{g(t)}\\
      =& p'(t)   \int_{\bf M} v^2 \ln v^2 d\mu_{g(t)}
         -4[p(t)-1]  \int_{\bf M}
            (|\nabla v |^2 +\frac{1}{4} R v^2) d\mu_{g(t)}\\
            &+ \{4[p(t)-1]- p^2(t)\}  \int_{\bf M} \frac{1}{4} R v^2
              d\mu_{g(t)}-\frac{3}{4} p^2(t)  \int_{\bf M}  R v^2 d\mu_{g(t)}.
 \eal
\]Notice that we have the following relations,
\[
\frac{4(p(t)-1)}{p'(t)} = \frac{4 t (T-t)}{T} \le T, \qquad  \frac{p^2(t)}{p'(t)}=T,
\]
\[
-T \le \frac{4(p(t)-1)-p^2(t)}{p'(t)} = \frac{4t(T-t)-T^2}{T} \le 0.
\]Hence we have
\[
\al & p^2(t)
\partial_t \ln \Vert u \Vert_{p(t)} \\
\le& p'(t) \bigg[  \int_{\bf M} v^2 \ln v^2  d\mu_{g(t)} -
\frac{4(p(t)-1)}{p'(t)} \int_{\bf M}
            (|\nabla v |^2 +\frac{1}{4} R v^2) d\mu_{g(t)}
+  T \sup R^-(x, t) \bigg]. \eal
\]Take $\epsilon$ such that
\[
\epsilon^2= \frac{4(p(t)-1)}{p'(t)} \le T \] in the
log-Sobolev inequality (\ref{reslogsob1}) (see [Section \ref{app},
appendix]), we deduce that
\[
p^2(t)
\partial_t \ln \Vert u \Vert_{p(t)} \le
 p'(t) \bigg[ -n \ln \sqrt{4 (p(t)-1)/p'(t)} + L(t) +  T \sup
 R^-(x, 0) \bigg],
 \]where  \[
 \al
 L(t)  &\doteq (t +\epsilon^2) \beta  + \alpha\\
&\le 2 T  \beta + \alpha \\
& \doteq L(T), \eal
 \] for some positive constants $\alpha=\alpha(A_0, B_0, \lambda_0, n)$ and
 $\beta=\beta(A_0, B_0, \lambda_0, n)$ from (\ref{reslogsob1}).
 Here we have used the fact that $\sup
 R^-(x, t) \le \sup
 R^-(x, 0)$.
Recall that $p'(t)/p^2(t)=1/T$ and $4(p(t)-1)/p'(t)=4t (T-t)/T$.
 Hence we have
 \[
\partial_t \ln \Vert u \Vert_{p(t)} \le
 \frac{1}{T} \bigg\{ - \frac{n}{2} \ln [4 t (T-t)/T] + L(T) + T \sup
 R^-(x, 0) \bigg\}.
 \]This yields, after integrating from $t=0$ to $t=T$, that
 \[
 \ln  \frac{\Vert u(\cdot, T) \Vert_\infty}{\Vert u(\cdot, 0)
 \Vert_1}
 \le - \frac{n}{2} \ln (4 T) + L(T) + T \sup
 R^-(x, 0)+n.
 \]Since
 \[
 u(x, T) = \int_{\bf M}  G(z, 0; x, T) u(z, 0) d\mu_{g(0)},
 \]
the above inequality implies that
 \begin{equation}
 \lab{ondiagbound}
 G(z, 0, x, T) \le \frac{\exp[L(T)+  T \sup
 R^-(x, 0)]}{(4  T)^{n/2}},
 \end{equation}
 where $L(T)$ is defined above as
 \[
 L(T) = 2  T \beta+\alpha.
\] As $T$ is arbitrary, we get the desired upper bound in  (\ref{odul}).  Note the constants  $\beta$ may have changed by a factor of $2$ and $\alpha$ has changed its value by $n$.

If $R(x, 0)>0$, then  it follows from the definition of
$\lambda_0$  in (\ref{lambda0}) that $\lambda_0>0$.  Moreover in (\ref{reslogsob1}),
we have $\beta=0$.
 So the above bound becomes
\begin{equation}
\label{ondiagbound2}
 G(z, 0; x, T) \le \frac{\exp(\alpha )}{(4 \pi T)^{n/2}},
 \end{equation} whence the upper bound in  (\ref{odulp}) follows.\\

Next we prove a lower bound. Let $t<t_0$  and $u=u(x, t) \equiv
G(x, t; x_0, t_0)$.   We claim that for
 a constant $C>0$,
 \[
 G(x_0, t; x_0, t_0) \ge  \frac{C}{\tau^{n/2}} e^{- \frac{1}{2 \sqrt{\tau}} \int^{t_0}_{t} \sqrt{t_0-s}
 R(x_0, s) ds} .
 \]where $\tau = t_0-t$ here and later in the proof.
  To prove this inequality,  define  a function $f$ by
 \[
 (4 \pi \tau)^{-n/2} e^{-f} = u.
 \]As a
 consequence of Perelman's differential Harnack inequality for the
 fundamental
 solution along any smooth space-time curve $\gamma (t)$ (see \cite[Corollary 9.4]{perelman1}), here we pick the curve
 $\gamma(t)$ to be the fixed point $x_0$, we have,
 \[
 -\partial_t f(x_0, t) \le \frac{1}{2} R(x_0, t) - \frac{1}{2 \tau} f(x_0, t).
 \]For any $t_2<t_1<t_0$, we can integrate the above inequality to get
 \[
 f(x_0, t_2) \sqrt{t_0-t_2} \le  f(x_0, t_1) \sqrt{t_0-t_1}  + \frac{1}{2} \int^{t_1}_{t_2} \sqrt{t_0-s}
 R(x_0, s) ds.
 \]When $t_1$ approaches $t_0$,  $f(x_0, t_1)$ stays bounded since
$G(x_0, t_1; x_0, t_0) (t_0-t_1)^{n/2} $ is bounded between two
positive constants, which is a direct consequence of the standard asymptotic formula for $G$
(for example, see \cite[Chapter 24]{chowetc3}).  Hence for any $t \le t_0$, we
have
 \[
 f(x_0, t)  \le  \frac{1}{2 \sqrt{t_0-t}} \int^{t_0}_{t} \sqrt{t_0-s}
 R(x_0, s) ds.
 \]Consequently
\be
\lab{Gx0x0}
G(x_0, t; x_0, t_0) \ge \frac{c}{(4 \pi \tau)^{n/2}} e^{-
\frac{1}{2 \sqrt{t_0-t}} \int^{t_0}_{t} \sqrt{t_0-s}
 R(x_0, s) ds} .
 \ee Notice that this bound is a global one that requires no curvature assumption. It also
holds on complete noncompact manifolds whenever maximum principle applies.
This will be important when we  study type I $\kappa$-solution in Section 4.

 We observe that $G(x_0, t; \cdot, \cdot)$ is a
solution to the standard heat equation coupled with Ricci flow,
which is the conjugate of the conjugate heat equation. i.e.,
 \[
 \Delta_z G(x,
t; z; l) - \partial_l G(x, t; z, l) =0,
\] here $\Delta_z$ is with respect to the metric $g(l)$. Therefore
it follows from \cite[Theorem 3.3]{z06} or \cite[Theorem
5.1]{ch2009} that, for $\delta>0, c_1, c_2>0$, and $y_0 \in \M$,
\[
G(x_0, t; x_0, t_0) \le c_1 G^{1/(1+\delta)} (x_0, t, y_0, t_0)
K^{\delta/(1+\delta)}  e^{c_2 d^2(x_0, y_0, t_0)/\tau},
\]where $K =\sup_{M \times [t/2, 0]} G(x_0, t, \cdot, \cdot)$.
By the on-diagonal upper bound
\[
K \le \frac{c B(\tau/2)}{\tau^{n/2}},
\]
this together with the on-diagonal lower  bound shows that, with $\delta = 1$,
\[
G(x_0, t; y_0, t_0) \ge c_1 \frac{B^{-1}(\tau/2)}{ \tau^{n/2}}
e^{-2 c_2 d(x_0, y_0, t_0)^2/\tau}  e^{- \frac{1}{ \sqrt{t_0-t}}
\int^{t_0}_{t} \sqrt{t_0-s}
 R(x_0, s) ds},
\] which is our desired bound.
\qed

%%%%%%%%%%%%%%%%%%%%%%%%%%%%%%%%%%%%%%%%%%%%%%%%%%%%%%%%%%%%%%%%%%%%%%%%%
\section{Full upper and lower bound with nonnegative Ricci curvature}
\label{sec:3}
%%%%%%%%%%%%%%%%%%%%%%%%%%%%%%%%%%%%%%%%%%%%%%%%%%%%%%%%%%%%%%%%%%%%%%%%%

In this section, we focus on the case of \textbf{nonnegative}
Ricci curvature. We establish Gaussian  upper and lower bounds for
the fundamental solution of the conjugate heat equation, these
bounds are in global nature. The result in this section can be
regarded as a generalization of Li and Yau's estimate (\cite{ly86}) from fixed metrics to metrics evolving under  the Ricci flow. Namely, we obtained sharp bounds for the fundamental solution of the conjugate heat equation on manifolds with nonnegative Ricci curvature. The main technique that we use here is much
different since now that the metric is evolving by the Ricci flow. The
essential new tools are an elliptic type Harnack inequality
for the fundamental solution and a uniform Sobolev inequality,
which was proven recently in \cite{ch2009, z06} and \cite{z07}.
Previous related results along the Ricci flow setting can be found
in \cite{guenther02, ni04, perelman1, z06}.

\medskip
Our main result of this section is the following,

\begin{theorem}\label{thm:3.1}
%th3.1
Let $({\bf M}, g(l))$, $l\in [0,T)$, be a solution the Ricci flow. Let
$G=G(x, l ; y, t)$ be the fundamental solution of the conjugate
heat equation. Assume that $({\bf M}, g(l))$ has nonnegative Ricci
curvature and is not Ricci-flat. Then for any $l, t \in [0, T)$,
$l <t$,  and $x, y \in {\bf M}$,  there exist positive constants  $c$ and
$c_n~($which only depends on $g_0$ and the dimension of $\bf M)$, and numerical
constants  $\eta_1$
, $\eta_2$
such that
\[
\frac{e^{-\eta_2 \Lambda_2(t)}}{c_n|B(x, \sqrt{t-l}, t)|_t}    e^{- d(x, y, t)^2/c(t-l)}   \le G(x, l; y, t) \le
 \frac{c_n e^{-\eta_1 \Lambda_1(t)} }{|B(x, \sqrt{t-l},
t)|_t}    e^{-c d(x, y, t)^2/(t-l)},
\]here $|B(x, r, t)|_t$ denote the volume of the ball $B(x, r, t)$ measured using $g(t)$, $\Lambda_1 (t)=\int_0^t \min R(\cdot,
s) ds$ and
 $\Lambda_2 (t)=\int_0^t \max R(\cdot,
s) ds$.
\end{theorem}

\vspace{.1in}

\proof Without loss of generality we assume that $l=0$. Our idea
is to first bound the fundamental solution $p=p(y, t; x, 0)$ of
the (forward) heat equation
\begin{equation}\label{he}
\partial_t u =\Delta_y u.
\end{equation}\\
Observe that
\[
p(y, t; x, 0) = G(x, 0; y, t),
\]where $G$ is the fundamental solution of the conjugate heat equation, hence a bound for $p$ is also a bound for $G$. We start with
Grigor'yan's method in
\cite{grigoryan97}.\\

\noindent {\bf Step 1.}  We first obtain monotonicity of certain
weighted $L^2$ norms of the solution. Let $u$ be a positive
solution to the equation (\ref{he}). Pick a weight function
$e^{\xi(y, t)}$ which we will specify later.  We compute that \be \lab{ddsu2}
 \frac{d}{d t} \int_{\M} u^2 e^\xi d\mu_{g(t)}(y) =
 \int_{\M} u^2 e^\xi \partial_t \xi d\mu_{g(t)}(y)+ \int_{\M} 2  u (\Delta u
 -\frac{R}{2} u) e^\xi d\mu_{g(t)}(y).
 \ee Note that
\[
\al
\int_{\M} u \Delta u e^\xi d\mu_{g(t)}(y) &= - \int_{\M} \nabla u \nabla
( u e^\xi) d\mu_{g(t)}(y)\\
&=- \int_{\M} \nabla u \nabla( u e^{\xi/2}  \   e^{\xi/2} )
d\mu_{g(t)}(y)
\\ &= - \int_{\M} \nabla u \left[ \nabla( u e^{\xi/2} )  \
e^{\xi/2} +  u e^{\xi/2}
\nabla  e^{\xi/2}  \right] d\mu_{g(t)}(y)\\
&= -\int_{\M} |\nabla( u e^{\xi/2} ) |^2 d\mu_{g(t)}(y) + \int_{\M} u^2
|\nabla e^{\xi/2}  |^2 d\mu_{g(t)}(y). \eal
\]Substituting this to the right hand side of (\ref{ddsu2}), we obtain
\[
\frac{d}{d t} \int_{\M} u^2 e^\xi d\mu_{g(t)}(y) \le  \int_{\M}
(\partial_t \xi + \frac{1}{2} |\nabla \xi|^2)  u^2 e^\xi d\mu_{g(t)}(y)-\int R u^2 e^\xi d\mu_{g(t)}(y) .
\]
If we choose $\xi$ such that
\[
\partial_t \xi + \frac{1}{2} | \nabla \xi|^2 \le 0,
\]then it follows that
\be
\lab{monou2exi}
 \int_{\M} u^2 e^\xi d\mu_{g(t_1)}(y) \bigg{|}_{t_1} \le
 e^{-(\Lambda_1(t_1)-\Lambda_1(t_2))}
 \int_{\M} u^2 e^\xi d\mu_{g(t_2)}(y) \bigg{|}_{t_2},
 \ee
 for $t_2 <t_1$.\\
%\medskip

\noindent {\bf Step 2.} With the above monotonicity
formula, we now use an idea from \cite{grigoryan97} to obtain
Gaussian upper bound for certain integral of $p(y,t; x,0)$. Fixing
a point $x \in \M$ and some positive constants $s$ and $r$, we
define
 \be
\lab{IRS} I_r(t) = \int_{{\M}-B(x, r, t)} u^2(y, t) d\mu_{g(t)}(y).
 \ee  We want to show that $I_r(t)$ has certain exponential decay for
 $u(y, t) = p(y, t; x, 0)$.  Take $A\ge 2$ and fix $t_0>0$, for $t<t_0$,
  we choose
 \[
 \al
 \xi =\xi(y, t)=
 \begin{cases}
  -\frac{(r-d(x, y, t))^2}{A (t_0-t)}, &\qquad d(x, y, t)
 \le r;\\
  0, &\qquad d(x, y, t) >r.
 \end{cases}
 \eal
 \]Then for $y \in B(x, r, t)$, we have
\[
\partial_t \xi + \frac{1}{2} | \nabla \xi|^2 =
- \frac{(r-d(x, y, t))^2}{A (t_0-t)^2} + \frac{2(r-d(x, y,
t))^2}{A^2 (t_0-t)^2} + \frac{2(r-d(x, y, t)) \partial_t d(x, y,
t)}{A (t_0-t)} \le 0,
\] here we used the nonnegativity of Ricci curvature and
hence that $\partial_t d(x, y,t) \le 0$.

By (\ref{monou2exi}), we have, for $t_2<t_1<t_0$,
\[
\int_{\M} u^2 e^\xi d\mu_{g(t)}(x) \bigg{|}_{t_1} \le
 \int_{\M} u^2 e^\xi d\mu_{g(t)}(x) \bigg{|}_{t_2} e^{-(\Lambda_1(t_1)-\Lambda_1(t_2))}.
 \]Since $\xi(y, t) =0$ when $d(x, y, t) \ge r$, this implies
 \[
 \al
 I_r(t_1) &= \int_{{\M}-B(x, r, t_1)} u^2(y, t_1) d\mu_{g(t_1)}(y) \le
 \int_{\M} u^2(y, t_1)  e^{\xi(y, t_1)} d\mu_{g(t_1)}(y)\\
 &\le
\int_{\M} u^2(y, t_2)  e^{\xi(y, t_2)} d\mu_{g(t_2)}(y) e^{-(\Lambda_1(t_1)-\Lambda_1(t_2))}.
\eal
\]For any number $\rho<r$, we can write this inequality as
\[
\al
 I_r(t_1) \le \left[  \int_{B(x, \rho, t_2)} u^2(y, t_2)  e^{\xi(y, t_2)} d\mu_{g(t_2)}(y)
   +
 \int_{{\M}-B(x, \rho, t_2)} u^2(y, t_2)
e^{\xi(y, t_2)} d\mu_{g(t_2)}(y) \right]  e^{-(\Lambda_1(t_1)-\Lambda_1(t_2))},
\eal
\]which shows that
\[
I_r(t_1) \le \left[ I_{\rho}(t_2) +
e^{-(r-\rho)^2/(A (t_0-t_2))} \ \int_{B(x, \rho, t_2)} u^2(y, t_2)  d\mu_{g(t_2)}(y) \right]
e^{-(\Lambda_1(t_1)-\Lambda_1(t_2))}.
\]Now we take $u(y, t) = p(y, t;
x, 0)$ be the fundamental solution.  For the fundamental solution, it holds that
\[
\al
&\int_{B(x, \rho, t_2)} u^2(y, t_2) d\mu_{g(t_2)}(y)
\\
&\le \int_{\M} p^2(y, t_2; x, 0) d\mu_{g(t_2)}(y) \le
 Q( t_2) \int_{\M} p(y, t_2; x, 0) d\mu_{g(t_2)}(y) \le   Q(t_2) e^{-\Lambda_1(t_2)},
 \eal
\]where $Q(t)$ is given by the righthand side of the on-diagonal bound in the previous
section. In particular, since  $R \ge 0$ and hence
$B_0=0$, we have $Q(t) = ct^{-n/2}$.  In the last inequality, we
have used
\[
\int_{\M} p(y, t_2; x, 0) d\mu_{g(t_2)}(y) \le e^{-\Lambda_1(t_2)},
\] due to the fact that $$\ppt \int_{\M} p(y, t; x, 0) d\mu_{g(t)}(y)=-
\int_{\M} Rp(y, t; x, 0) d\mu_{g(t)}(y) \le -\min R(\cdot, t) \int_{\M} p(y, t; x, 0) d\mu_{g(t)}(y) .$$

Thus we reach the following inequality
\[
I_r(t_1) \le \left[ I_{\rho}(t_2) + e^{-(r-\rho)^2/(A(t_0-t_2))} \
Q(t_2)  e^{-\Lambda_1(t_2)}\right]
e^{-(\Lambda_1(t_1)-\Lambda_1(t_2))}.
\]Observe the above definition of $I_r(t)$ is independent of $t_0$ and
$\xi$. So we can take $t_0=t_1$ and it follows that
\be
\lab{IrIrho}
 I_r(t_1) \le
 I_{\rho}(t_2)  e^{-(\Lambda_1(t_1)-\Lambda_1(t_2))} +
 e^{-(r-\rho)^2/(A(t_1-t_2))} \ Q(t_2) e^{-\Lambda_1(t_1)},
\ee where $r>\rho$, $t_1>t_2$ and $A \ge 2$.\\

Now fixing $r, t>0$,  we define two
sequences $\{r_k\}$ and $\{t_k\}$ as in \cite{grigoryan97},
\[
r_k = \left( \frac{1}{2} + \frac{1}{k+2} \right) r, \qquad t_k =
\frac{t}{a^k}, \quad k=0, 1, 2, ...
\]where $a>1$ will be chosen later (notice that $t_0$, $t_1$ and $t_2$
defined here are not related to those appeared above). Applying
(\ref{IrIrho}), we deduce \be \lab{Irk}
 I_{r_k}(t_k) \le
I_{r_{k+1}}(t_{k+1}) e^{-(\Lambda_1(t_k)-\Lambda_1(t_{k+1}))} +
e^{-(r_k-r_{k+1})^2/(A(t_k-t_{k+1}))} \ Q( t_{k+1})
e^{-\Lambda_1(t_k)}.
 \ee

Remember that $I_{r_k}(t_k) = \int_{{\M}-B(x, r_k, t_k)} p^2(y, t_k; x,
0) d\mu_{g(t)}(y)$. When $k \to \infty$, $t_k \to 0$ and $p(y, t_k; x, 0)  \to
\delta(y, x)$ which is concentrated at the point $x$. Hence $\lim_{k
\to \infty} I_{r_k}(t_k) =0.$  This argument can easily be made
rigorous by approximating $p$ with  regular solutions whose initial
value is supported in $B(x, r/2, 0)$.

After applying iterations to (\ref{Irk}), we obtain that
\[
I_r(t) =  I_{r_0}(t_0) \le e^{-\Lambda_1(t)}   \sum^\infty_{k=0}
Q( t_{k+1}) e^{-(r_k-r_{k+1})^2/(A(t_k-t_{k+1}))} .
\]Using the relation
\[
r_k - r_{k+1} \ge r/(k+2)^2, \qquad t_k-t_{k+1} = (a-1) t/a^{k+1},
\]we arrive at
\[
I_r(t)   \le e^{-\Lambda_1(t)} \sum^\infty_{k=0} Q( t_{k+1}) \exp
\left(-\frac{ a^{k+1}  \ r^2}{ (k+3)^4 \ (a-1) \ A t} \right).
\]
Substituting $Q(t_{k+1}) =c \frac{a^{(k+1)n/2} }{t^{n/2}}$ into
the last inequality concerning $I_r(t)$, we deduce
\[
I_r(t)  \le e^{-\Lambda_1(t)} \frac{ c}{t^{n/2}} \sum^\infty_{k=0}
 a^{(k+1)n/2}  \exp
\left(-\frac{ a^{k+1}  \ r^2}{ (k+3)^4 \ (a-1) \ A t} \right).
\]By making the constant $a$ sufficiently large and taking
$r^2\ge \frac14t$,
 it leads to that
\be
 \lab{IRS2}
 I_r(t) = \int_{{\M}-B(x, r, t)} p^2(y, t, x, 0) d\mu_{g(t)}(y)
 \le e^{-\Lambda_1(t)} \frac{ c}{t^{n/2}}
e^{-c_1 r^2/t} \ee for some positive constants $c $ and $c_1$.

Let $x_0, y_0 \in {\M}$ be two points such that $d(x_0, y_0, t) \ge
\sqrt{t}$.
 Then
\[
B(y_0, \sqrt{t/4}, t) \subset {\M}-B(x_0, r, t),
\] where $r=d(x_0, y_0, t)/2$.  Hence, it follows from  (\ref{IRS2}) that,
there exists $z_0 \in B(y_0, \sqrt{t/4}, t)$ such that
\[
p^2(z_0, t, x_0, 0) |B(y_0, \sqrt{t/4}, t)|_t \le \frac{
c}{t^{n/2}} e^{-c_1 r^2/t}  e^{-\Lambda_1(t)},
\]
i.e.,
\[
p^2(z_0, t; x_0, 0) \le \frac{c e^{-c_1 d(x_0, y_0, t)^2/t}}{
|B(x_0, \sqrt{t/4}, t)|_{t}  t^{n/2}}  e^{-\Lambda_1(t)}.
\]By the classical volume comparison theorem, this implies that
\be \label{eq:pptdec} p^2(z_0, t; x_0, 0) \le \frac{c_n e^{-c_2
d(x_0, y_0, t)^2/t}}{ |B(x_0, \sqrt{t}, t)|^2_t }
e^{-\Lambda_1(t)}, \ee for some positive constants $c_n$ and
$c_2$.
\\

\noindent {\bf Step 3.} Next, let us recall that $p(z, t; x_0, 0)$
is a solution to the heat equation. i.e.
 \[ \Delta_z p(z,
t; x_0; 0) - \partial_t  p(z, t; x_0, 0) =0.
\] Using \cite[Theorem 3.3]{z06}
(also see \cite{ch2009}), for any $\delta>0$, there exist positive
constants $c_3$ and  $c_4$, \be \lab{pMe} p(y_0, t; x_0, 0) \le
c_3 p^{1/(1+\delta)} (z_0, t; x_0, 0) K^{\delta/(1+\delta)} e^{
c_4 d^2(z_0, y_0, t)/t}, \ee where $K =\sup_{{\M} \times [t/2, t]}
p(\cdot, \cdot; x_0, 0)$. By Theorem \ref{thm:2.1}, there exists a
constant $c>0$, such that
\[
K \le  \frac{c}{t^{n/2}}.
\]It follows from (\ref{eq:pptdec}), (\ref{pMe}) and volume
comparison theorem that
\[
p(y_0, t; x_0, 0)^2 \le \frac{c_n e^{-c d(x_0, y_0,
t)^2/t}}{|B(x_0, \sqrt{t}, t)|^2_t } e^{-\Lambda_1(t)}.
\]
This shows, since $p(y_0, t; x_0, 0) = G(x_0, 0; y_0, t)$,
\[
G(x_0, 0; y_0, t) \le \frac{c_n e^{-c d(x_0, y_0,
t)^2/{t}}}{|B(x_0, \sqrt{t}, t)|_t } e^{-\Lambda_1(t)/2}.
\]
Since $x_0$ and $y_0$ are arbitrary, this proves the desired upper
bound.
\\

\noindent {\bf Step 4.} Next we show that a lower bound follows
from the upper bound. Recall the notation $\Lambda_2 (t)=\int_0^t \max R(\cdot,
s) ds$ and notice that
\[
\al
\frac{d}{d t} \int p(x, t; x_0, 0) d\mu_{g(t)}(x)d\mu_{g(-s)}(y) =& - \int R(x, t)  p(x,
t; x_0, 0) d\mu_{g(t)}(x) \\
\ge& - \max R(\cdot, t)   \int  p(x, t; x_0, 0)
d\mu_{g(t)}(x).
\eal
\]Hence
\[
\int p(x, t; x_0, 0) d\mu_{g(t)}(x) \ge e^{- \Lambda_2 (t) }.
\]

For $\beta>0$ that we will fix later, the upper bound implies
\[
\al
&\int_{B(x_0, \sqrt{ \beta t}, t)} p^2(x, t; x_0, 0) d\mu_{g(t)}(x)\\
&\ge \frac{1}{|B(x_0, \sqrt{ \beta t}, t)|_t}
\left(  \int_{B(x_0, \sqrt{ \beta t}, t)} p(x, t; x_0, 0) d\mu_{g(t)}(x) \right) ^2\\
&=\frac{1}{|B(x_0, \sqrt{ \beta t}, t)|_t}
\left( e^{-\Lambda_2 (t)}- \int_{\M-B(x_0, \sqrt{ \beta t}, t)} p(x, t; x_0, 0) d\mu_{g(t)}(x) \right) ^2\\
&\ge \frac{1}{|B(x_0, \sqrt{ \beta t}, t)|_t} \left( e^{-\Lambda_2
(t)}- \int_{\M-B(x_0, \sqrt{ \beta t}, t)} \frac{c_n}{|B(x_0,
\sqrt{t}, t)|_t}
  \   e^{-c \ d(x_0, x, t)^2/t} d\mu_{g(t)}(x) \right) ^2.
\eal
\]Since the Ricci curvature is nonnegative,
one can use the volume doubling property to compute that, for
$\beta=2 (\Lambda_2 (t)+C)/c$ with $C$ large enough, we have
\[
\al& \int_{\M-B(x_0, \sqrt{ \beta |t|}, t)} \frac{c_n}{|B(x_0,
\sqrt{t}, t)|_t}
  \   e^{-c \ d(x_0, x, t)^2/t} d\mu_{g(t)}(x) \\
  & \le \int_{\M-B(x_0, \sqrt{ \beta t}, t)} \frac{c_n}{ 2|B(x_0,
\sqrt{t}, t)|_t}
  \   e^{- c \ d(x_0, x, t)^2/(2t)} d\mu_{g(t)}(x) e^{-c \beta/2} \le \frac12 e^{-c \beta/2}.
  \eal
\]Hence there exists $x_1 \in B(x_0, \sqrt{ \beta |t|}, t)$ such that
\[
p(x_1, t; x_0, 0) \ge \frac{1}{|B(x_0, \sqrt{\beta  t}, t)|_t}
(e^{- \Lambda_2 (t)}-\frac12e^{-c \beta/2}).
\]Recall $ C+\Lambda_2 (t) = c \beta/2$, we deduce
\[
p(x_1, t; x_0, 0) \ge \frac{1}{2 |B(x_0, \sqrt{\beta  t}, t)|_t}
e^{- \Lambda_2 (t)}.
\]Using volume comparison theorem again, we have
\[
|B(x_0, \sqrt{\beta  t}, t)|_t \le  \beta^{n/2} |B(x_0, \sqrt{ t},
t)|_t\le (C +\Lambda_2 (t))^{n/2} |B(x_0, \sqrt{ t}, t)|_t,
\]where the constant $C$ may have changed its value. Thus
\[
p(x_1, t; x_0, 0) \ge \frac{1}{2(C+\Lambda_2 (t))^{n/2} |B(x_0,
\sqrt{t}, t)|_t} e^{- \Lambda_2 (t) }\ge \frac{\eta_1}{|B(x_0,
\sqrt{t}, t)|_t}  e^{- (1+\eta_2)\Lambda_2 (t) }.
\]

Now we can use the same approach as we derive the upper bound. As in (\ref{pMe}), for any $y_0 \in {\bf M}$, we have  \be
\lab{pMe2} p(x_1, t; x_0, 0) \le C p^{1/(1+\delta)} (y_0, t;  x_0,
0) K^{\delta/(1+\delta)} e^{ c d^2(x_1, y_0, t)/t}. \ee  Since
$d^2(x_1, y_0, t)/t \le 2\beta+c_2 d^2(x_0, y_0, t)/t$, by taking
$\delta=1$ and using the volume comparison theorem, we have
\[
G(x_0, 0; y_0, t)=p(y_0, t;  x_0, 0) \ge \frac{C e^{- c  d^2(x_0,
y_0, t)/t}}{|B(x_0, \sqrt{t}, t)|_t}e^{-\eta_3\Lambda_2 (t) }.
\]
This is a lower bound which matches the upper bound except for
constants. \qed

\bigskip

 \section{Applications to Type I $\kappa$-solutions}
 \label{sec:4}

In this section, as an application of our previous bounds on the
fundamental solution, we shall obtain a classification for
backward limit of $\kappa$-solutions without assuming nonnegative
curvature operator. In dimension $3$, as a result of Hamilton-Ivey
curvature pinching estimate, $\kappa$-solutions have nonnegative
curvature operator (hence nonnegative sectional curvature).
Perelman \cite{perelman1} gave a classification for all such
solutions. For dimension at least $4$, a priori, we can no longer
assume nonnegative curvature operator since there is no
Hamilton-Ivey type estimate. Moreover the $\kappa$-solutions maybe
complete non-compact, so we need to prove the bounds for the
fundamental solution in the complete non-compact setting. Our
method here is similar to the one we used in previous two sections
where we have a closed solution of the Ricci flow.\\

For convenience and without loss of generality, we take the final
time $T_0=1$ for the ancient solution throughout the section; we
also take $D_0 \ge 1$. The conjugate heat equation is \be
\lab{eqconj}
 \Delta u - R u - \partial_\tau u=0,
\ee     here $\tau = - t$,  $\Delta$ and $R$ are the
Laplace-Beltrami operator and the scalar curvature with respect to
$g(t)$. This equation, coupled with the initial value $u_{\tau=0}
= u_0$ is well posed if $\M$ is compact or if curvature is
bounded and $u_0$ is bounded \cite{guenther02}.

We use $G=G(x, \tau; x_0, \tau_0)$ to denote the heat kernel
(fundamental solution) of (\ref{eqconj}),  here $\tau>\tau_0$ and
$x, x_0 \in \M$. In the rest of this paper, if $\bf M$ is complete
non-compact, $G$ is meant to be the minimal fundamental solution.
Our main technical result of the section is the following,

\begin{lemma}
\label{le1.1} Let $(\M \sp n, g(t))$, $t \in (-\infty, 0]$, be a
$\kappa$-solution to the Ricci flow. Then there exist positive
numbers $a_1$ and $b_1$ which only depend on $n$, $\kappa$ and
$D_0$, such that for all $x, x_0 \in \M$,  and $\tau = - t
>0$, we have
\[
\al & G(x, \tau; x_0, 0) \le \frac{a_1}{\tau^{n/2}},\\
& G(x_0, \tau; x_0, 0) \ge \frac{1}{a_1 \tau^{n/2}}, \eal
\] here $G(x, \tau; x_0, 0)=G(x, t; x_0, 0)$ as in Theorem
\ref{thm:2.1}.
\end{lemma}

\proof The proof of the lemma is similar to that in Section 2.
Comparing with that case, we have two new ingredients coming from
type I $\kappa$-solutions. The first one is the non-collapsing
condition on all scales. The second one is the bound on
curvatures. These new ingredients allow us to obtain a better
estimate. It is convenient to work with the reversed time $\tau$.
Notice that the Ricci flow becomes a backward flow with respect to
$\tau$ and the conjugate heat equation is now forward conjugate
heat equation.\\

{\bf Step 1.} Since $Ric(x, t) \ge -\frac{D_0}{1+|t|}$ (this $D_0$
differs by a dimensional constant with the one in Definition 1.2),
it is well known (e.g. \cite[Theorem 3.1]{sc92}) that the following
Sobolev inequality holds: Let $B(x, r, t)$ be a proper sub-domain
for $(\M, g(t))$, for all $v \in W^{1, 2}_0(B(x, r, t))$, there
exists $c_1, c_2$ depending only on dimension $n$ such that,
\be
 \lab{SobRic>0}
  \bigg{(}  \int v^{2n/(n-2)} d\mu_{g(t)}
\bigg{)}^{(n-2)/n} \le \frac{c_1  r^2 e^{c_2 r \sqrt{
D_0/|t|}}}{|B(x, r, t)|^{2/n}_t}  \int \left[ | \nabla v |^2  +
r^{-2} v^2 \right] d\mu_{g(t)}. \ee

In this section, we always take $r = c \sqrt{|t|}$, for some
$c<1$. By the assumption that $R(x, t) \le \frac{D_0}{1+|t|}$,
$D_0 \ge 1$, and the $\kappa$-non-collapsing property, we have
\[
| B(x, \sqrt{|t|},  t)|_t \ge \kappa D^{-n/2}_0 |t|^{n/2}.
\]Therefore the above Sobolev inequality becomes
\be \lab{SobT} \bigg{(}  \int v^{2n/(n-2)} d\mu_{g(t)}
\bigg{)}^{(n-2)/n} \le \frac{c_3(D_0, n)}{\kappa^{2/n}}  \int
\left[ | \nabla v |^2  + |t|^{-1} v^2 \right] d\mu_{g(t)}, \ee for
all $v \in W^{1, 2}_0(B(x, \sqrt{|t|},  t))$.\\

Next we show that, under the assumptions of the theorem, $({\M},
g(t))$ possess a space-time doubling property, i.e., the distance
between two points at different times $t_1$ and $t_2$ are
comparable if $t_1$ and $t_2$ are comparable. For any $x_1, x_2
\in \M$, let $\bf \gamma$ be the shortest geodesic connecting the
two. Then
\[
- \sup_{\bf \gamma} \int_{\bf \gamma} Ric (T, T) ds \le
\partial_t d(x_1, x_2, t) \le - \inf_{\bf \gamma} \int_{\bf \gamma} Ric
(T, T) ds,
\] here $T$ is the unit tangent vector of $\bf \gamma$. By curvature
assumption, it holds
\[
| Ric(x, t) | \le \frac{c D_0}{1+|t|}.
\]Therefore
\[
 -\frac{c D_0}{1+|t|} d(x_1, x_2, t) \le \partial_t d(x_1, x_2, t) \le \frac{c D_0}{1+|t|} d(x_1, x_2, t).
\]After integration, we arrive at
\be \lab{distdouble} \left( |t_1|/|t_2| \right)^{c D_0} \le d(x_1,
x_2, t_1)/d(x_1, x_2, t_2) \le \left( |t_1|/|t_2| \right)^{-c D_0},
\ee  for all $t_2<t_1<0$. Note that the above inequality is of
local nature. If the distance is not smooth, one can just
shift one point, say $x_1$,  slightly and then obtain the same
integral inequality by taking limits.\\

Similarly, for any $x\in {\bf M}$ and fixed $t_1$, we have
\[
\bigg|\partial_t \int_{B(x, \sqrt{|t_1|}, t_1)} d\mu_{g(t)}\bigg|
= \bigg|- \int_{B(x, \sqrt{|t_1|}, t_1)} R(y, t) d\mu_{g(t)}\bigg|
\le \frac{c D_0}{1+|t|}
 \int_{B(x, \sqrt{|t_1|}, t_1)} d\mu_{g(t)}.
\]Upon integration, we know that the volume of the balls
\be
\lab{voldouble}
 |B(x,  \sqrt{|t_3|}, t_4)|_{t} \equiv
Vol_{g(t)} \{ y \, | \, d(x, y, t_4) < \sqrt{|t_3|} \}
 \ee
are all comparable for $t_3, t_4, t \in [t_2, t_1]$, provided that
$t_1$ and $t_2$ are comparable.

Let $u$ be a positive solution to (\ref{eqconj}) in the region
\[
Q_{\sigma r}(x, \tau) \equiv \{ (y, s) \ | \ y \in {\bf M}, \tau-(\sigma
r)^2 \le s \le \tau, \ d(y, x, -s) \le \sigma r \},
\]here $r=\sqrt{|t|}/8>0$, $2 \ge \sigma \ge 1$.
Given any $p \ge 1$, it is clear that
\be
\lab{equp}
\Delta u^p - p R u^p- \partial_\tau u^p \ge 0.
\ee

Let $\phi: [0, \infty) \to [0, 1]$ be a smooth function such that
$|\phi'| \le 2/((\sigma-1) r)$, $\phi' \le 0$, $\phi(\rho) =1$
when $0 \le \rho \le r$, and $\phi(\rho)=0$ when $\rho \ge \sigma r$.
Let $\eta: [0, \infty) \to [0, 1]$ be a smooth function such that
$|\eta'| \le 2/((\sigma-1) r)^2$, $\eta' \ge 0$, $\eta \ge 0$,
$\eta(s) =1$ when $\tau-r^2 \le s \le \tau$, and $\eta(s)=0$ when $s
\le \tau- (\sigma r)^2$.  Define a cut-off function $\psi =
\phi(d(x, y, -s)) \eta(s)$.\\

Let $w= u^p$ and using $w \psi^2$ as a test function for (\ref{equp}),
we deduce
\be
\lab{dwdw}
\int \nabla (w \psi^2) \nabla w d\mu_{g(-s)}(y)ds + p \int R w^2 \psi^2
d\mu_{g(-s)}(y)ds \le - \int (\partial_s w) w \psi^2  d\mu_{g(-s)}(y)ds.
\ee  Direct calculation yields
\[
\int \nabla (w \psi^2) \nabla w d\mu_{g(-s)}(y)ds = \int  | \nabla (w
\psi)|^2 d\mu_{g(-s)}(y)ds - \int  |\nabla \psi|^2 w^2 d\mu_{g(-s)}(y)ds.
\]Next we estimate the righthand side of (\ref{dwdw}),
\[
\al
&- \int (\partial_s w) w \psi^2  d\mu_{g(-s)}(y)ds\\
 = &\int w^2 \psi
\partial_s \psi d\mu_{g(-s)}(y)ds + \frac 1 2 \int (w \psi)^2 R d\mu_{g(-s)}(y)ds
- \frac{1}{2} \int (w \psi)^2 d\mu_{g(-\tau)}(y).
\eal
\]Observe that
\[
 \partial_s \psi = \eta(s)   \phi'(d(y, x, -s)) \partial_s d(y, x, -s)
+ \phi(d(y, x, -s)) \eta'(s).
\]Note also $|\phi'| \le 2/(1-\sigma) r$ and
\[
|\partial_s d(y, x, -s)| \le |\int_{\bf \gamma} Ric(T, T) dl| \le
\frac{Cr}{s} \le \frac{C}{r},
\]here ${\bf \gamma} = {\bf \gamma}(l)$ is a minimum geodesic connecting $y$ and $x$, parameterized by
arc length.  We have also used the fact that $ r^2 \le c s$  above,
which holds since $s \in [|t|/2, |t|]$ by our choice. Therefore we have
\[
| \partial_s \psi| \le C/[(1-\sigma) r]^2.
\]
Hence it follows that
\be
\lab{dwdw2}
\al
 &- \int (\partial_s w) w \psi^2  d\mu_{g(-s)}(y)ds \\
 \le &\frac{C}{[(1-\sigma) r]^2} \int w^2 \psi d\mu_{g(-s)}(y)ds +
 \frac 1 2 \int (w \psi)^2 R d\mu_{g(-s)}(y)ds
- \frac{1}{2} \int (w \psi)^2 d\mu_{g(-\tau)}(y). \eal \ee
Combining (\ref{dwdw}) and (\ref{dwdw2}), we obtain, in view of
$p \ge 1$ and $|R| \le \frac{C}{r^2}$,
\be \lab{dwdw3} \int  |
\nabla (w \psi)|^2 d\mu_{g(-s)}(y)ds + \frac{1}{2} \int (w \psi)^2
d\mu_{g(-\tau)}(y)d\mu_{g(-s)}(z) \le \frac{c(1+p)}{(\sigma-1)^2
r^2} \int_{Q_{\sigma r (x, \tau)}} w^2 d\mu_{g(-s)}(y)ds.
\ee   By
H\"older's inequality, we have
\be \lab{2/n} \int (\psi w)^{2(1+(2/n)}
d\mu_{g(-s)}(y) \le \bigg{(}   \int (\psi w)^{2n/(n-2))}
d\mu_{g(-s)}(y) \bigg{)}^{(n-2)/n} \bigg{(} \int (\psi w)^2
d\mu_{g(-s)}(y) \bigg{)}^{2/n}. \ee
\medskip

We claim that the diameter $d$ of ${\M}$ at time $t$ is a least $c
\sqrt{|t|}$ for some $c=c(n, \kappa, D_0)>0$. Without loss of
generality, we can assume $d \le \sqrt{|t|}$. Since $ Ric(x, t)
\ge -\frac{D_0}{|t|}$, it follows from the classical volume
comparison theorem and $\kappa$ non-collapsed assumption that
\[
\kappa D^{-n/2}_0 |t|^{n/2} \le |B(x, \sqrt{|t|}, t)|_{t} = |B(x,
d, t)|_{t} \le e^{ c \sqrt{D_0}} d^n.
\] Hence
the claim follows.\\

By the distance doubling property (\ref{distdouble}),   $B(x,
\sigma r, -s)$ is a proper sub-domain of ${\bf M}$, $s \in
[\tau-(\sigma r)^2, \tau]$, where  $r=\sqrt{|t|}/C$ for some
sufficiently large number $C$, for simplicity, we just take $C=8$.
 By the Sobolev inequality (\ref{SobT}), it holds
\[
\bigg{(}  \int (\psi w)^{2n/(n-2)} d\mu_{g(-s)}(y) \bigg{)}^{(n-2)/n} \le
c(\kappa, D_0)
\int [| \nabla (\psi w) |^2  + r^{-2} (\psi w)^2] d\mu_{g(-s)}(y),
\]for $s \in [t-(\sigma r)^2, t]$. Substituting  this and (\ref{dwdw3}) to (\ref{2/n}),
 we arrive at the estimate
\[
\int_{Q_r(x, \tau)} w^{2 \theta} d\mu_{g(-s)}(y) ds \le c(\kappa,
D_0) \bigg{(} \frac{1+p}{(\sigma-1)^2 r^2} \int_{Q_{\sigma r}(x,
\tau)} w^2 d\mu_{g(-s)}(y) ds \bigg{)}^{\theta},
\]with $\theta = 1+(2/n)$. Now we apply the above inequality repeatedly with
the parameters $\sigma_0=2, \sigma_i=2- \Sigma^i_{j=1} 2^{-j}$ and
$p=\theta^i$. This shows a $L^2$ mean value inequality \be
\lab{L2mvi} \sup_{Q_{r/2}(x, \tau)} u^2 \le \frac{c(\kappa,
D_0)}{r^{n+2} } \int_{Q_r(x, \tau)} u^2 d\mu_{g(-s)}(y). \ee

This inequality also holds if one replaces $r$ by any positive
number $r_1 <r\le \sqrt{|t|}$ since $|B(x, r_1, t)| \ge k c_n
|B(x, r, t)| (\frac{r_1}{r})^n \ge c r_1^n$ by the volume
comparison theorem.  Then one can derive (\ref{L2mvi}) using Moser iteration.\\

  Now we use a generic trick of P. Li and R. Schoen (\cite{ls84}),
  recall that they only use the doubling property of metric balls, we arrive at the following $L^1$ mean value inequality
\[
\sup_{Q_{r/2}(x, \tau)} u \le  \frac{c(\kappa, D_0)}{r^{n+2} }
\int_{Q_r(x, \tau)} u d\mu_{g(-s)}(z)ds.
\]We further remark that the doubling constant is uniform based on our discussion in
last paragraph.\\

Now we take  $u(x, \tau) = G(x, \tau; x_0, 0)$. Note that
$\int_{\bf M} u(z, s) d\mu_{g(-s)}(z) =1$ and $r=\sqrt{|t|}$, we have \be
\lab{Gdbound} G(x, \tau; x_0, 0) \le \frac{c(\kappa,
D_0)}{|t|^{n/2}}. \ee This proves the upper bound.\\

The lower bound is a consequence of (\ref{Gx0x0}) in Theorem
\ref{thm:2.1} and the curvature bound $|R(x, t)| \le
\frac{D_0}{|t|}$. Recall that there is no curvature condition
required in the proof of (\ref{Gx0x0}) for compact manifolds,
hence it carries over to complete non-compact manifolds with
bounded curvature where the maximum principle holds. As before,
set $\tau=-t$, we have
\[
G(x_0, \tau; x_0, 0) \ge \frac{c}{(4 \pi \tau)^{n/2}} e^{-
\frac{1}{2 \sqrt{-t}} \int^{0}_{t} \sqrt{-l}
 R(x_0, l) dl} .
\] Hence the desired lower bound follows from $|R(x, t)|
\le \frac{D_0}{|t|}$ after integration. \qed

Now we state our main theorem of this section.

\begin{theorem}
\lab{thbacklimit} Let $({\bf M}, g(t))$, $t \in (-\infty, 0]$, be
a non-flat, type I $\kappa$-solution to the Ricci flow for some
$\kappa>0$. Then there exist a sequence of points $\{q_k \}
\subset \M$, a sequence of times $t_k \to -\infty$, $k=1, 2, ...$,
and a sequence of re-scaled metrics
\[
g_k(x, s) \equiv |t_k|^{-1}  g(x, t_k + s |t_k| )
\]
%$g_k(x, s) \equiv |t_k|^{-1} g(x, t_k + s |t_k|)$
around $q_k$, such that $({\bf M}, g_k, q_k)$ converge to a
non-flat gradient shrinking Ricci soliton in $C^\infty_{loc}$
topology.
\end{theorem}

\proof
 In the special case that $\M$ is a dilation limit of a Type
 I maximal solution, A. Naber \cite{naber07}
   proved that $\M$ is a gradient shrinking Ricci
 soliton. However it is not clear if $\M$ is non-flat.
 Our proof here works for both compact and complete noncompact cases.\\

  By the $\kappa$-non-collapsed assumption and curvature bound $|Rm(\cdot,
  t)|
 \le \frac{D_0}{1+|t|}$, it follows from Hamilton's compactness theorem
 (see \cite{Hcomp}) that there exist a sequence of
 time $t_k$, with $\tau_k \equiv |t_k| \to \infty$, such that
 the following statement holds:

 for any fixed point $x_0 \in \M$, the pointed manifolds $({\M},
g_k, x_0)$ with metrics
 \[
 g_k \equiv \tau^{-1}_k g(\cdot, - s \tau_k)
 \] converge
  to a pointed manifold $({\M}_\infty, g_\infty(\cdot, s), x_\infty)$
  in $C^\infty_{loc}$ topology, here $s>0$.\\

 We first prove that $g_\infty$ is a gradient shrinking Ricci soliton.
 For $x \in \M$ and $s \ge 1$, let
 \[
 u_k=u_k(x, s) \equiv \tau^{n/2}_k \ G(x, s \tau_k; x_0, 0),
 \]here $G$ is the fundamental solution of the conjugate heat equation and $x_0$ is a fixed point.
 We choose $q_k=x_0$ for the above re-scaled manifolds.  By Lemma \ref{le1.1},
 we know that there exist a uniform positive constant $U_0$, such that $u_k(x, s) \le U_0$ for all $k=1, 2, ...$ and
 $x \in \M$.   Note that $u_k$ is a
 positive solution of the conjugate heat equation  on
 $({\M}, g_k(s))$, i.e.,
 \[
 \Delta_{g_k} u_k - R_{g_k} u_k - \partial_s u_k  =0.
 \]For any compact time interval in $(0, \infty)$,
 $u_k$ are uniformly bounded, moreover, $R_{g_k}$ and $Rm_{g_k}$ are uniformly
 bounded. It follows from the standard parabolic theory
that $u_k$ is H\"older continuous uniformly with respect to $g_k$.
Hence there exists a subsequence, still denoted as $\{
 u_k \}$, which converge (in $C^\alpha_{loc}$ topology)
 to a $C^\alpha_{loc}$ function $u_{\infty}$ on $({\M}_\infty,
g_\infty(s), y_\infty)$.\\

 It is easy to see that $u_\infty$ is a weak solution of
 the conjugate heat equation on  $({\M}_\infty, g_\infty(s))$, i.e.
 \[
 \int \int \left( u_\infty \Delta \phi - R_\infty u_\infty \phi + u_\infty \partial_s \phi \right)
 d\mu_{g_\infty(s)}ds =0,
 \]for all $\phi \in C^\infty_0({\M}_\infty \times (-\infty, 0])$.

 Again by standard parabolic theory, the function $u_\infty$, being bounded on compact
 time intervals,
 is a smooth solution of the conjugate heat equation on $({\M}_\infty, g_\infty(s), y_\infty)$.
 We only need to show that $u_\infty$ is not zero.

 By Lemma \ref{le1.1}, there exists
 a constant $a_1>0$, for all $\tau \ge 1$,
 \[
 G(x, \tau; x_0, 0) \ge \frac{1}{a_1 \tau^{n/2}}.
 \]
Since $u_k$ is a dilation of $G(x, \tau; x_0, 0)$,
  we derive that $u_k(x_0, s) \ge \frac{1}{2^na_1}>0$ for $s \in [1, 4]$.
   Hence $u_\infty(x_0, s) \ge \frac{1}{2^na_1}>0$.
Then the maximum principle yields that $u_\infty$ is positive everywhere.

 Let us recall that, for each $u_k$, Perelman's $W$-entropy  is defined as
 \[
 W_k(s) = W(g_k, u_k, s) =
 \int \left[ s ( |\nabla f_k|^2 + R_k) +f_k - n \right] u_k d\mu_{g_k(s)},
 \]where $f_k$ is determined by the relation
 \[
 (4 \pi s)^{-n/2} e^{-f_k} = u_k;
 \]and $R_k$ is the scalar curvature with respect to $g_k$.
 By the uniform upper bound for $u_k$, we know that there exist $c_0>0$ such that
 \be
 \lab{fk>c}
 f_k = - \ln u_k - \frac{n}{2} \ln ( 4 \pi s) \ge - c_0,
 \ee  for all $k=1, 2, ...$ and $s \in [1, 4]$.  Here the choice of
 the time
 interval of $s$
 is just for convenience, in fact any finite time interval also
 works.\\

 If $\M$ is noncompact, one needs to justify that the integral
 $W_k(s)$ is finite.
 For fixed $k$, $u_k$ has a generic Gaussian upper and lower bound
 with coefficients
 depending on $\tau_k$, curvature tensor and their derivatives,
 as shown in \cite{guenther02}.
  Since the curvatures are all bounded, the term $f_k u_k$ (which is
  essentially
 $-u_k \ln u_k$) is integrable. By
 \cite[(26.94)]{chowetc3}, at each fixed time level,
 $|\nabla f_k|^2 u_k =
 \frac{|\nabla u_k|^2}{u_k}$
 is also integrable.
  Hence $W_k(s)$ is well defined.\\

 Since $\int_{\M} u_k dg_k =1$, by (\ref{fk>c}) and
 $R_k(\cdot, s) \ge - c/s$, we know that there exists $c_1>0$,
 such that
 \be
 \lab{wklowerb}
 W_k(s) \ge -c_1,
 \ee  for all $k=1, 2...$ and $s \in [1, 4]$.

 Recall that $W$ is invariant under proper scaling, i.e.,
 \[
 W_k(s) =W(g_k, u_k, s)=W(\tau_k g_k, u, s \tau_k)= W(g, u, s \tau_k),
 \]where $u=u(x, l) =G(x, l, x_0, 0)$. According to \cite{perelman1},
 \be
 \lab{dwds}
 \frac{d W_k(s)}{ds} =
 -  2 s \int | Ric_{g_k} + Hess_{g_k} f_k - \frac{1}{2s}
 g_k |^2 u_k d\mu_{g_k(s)} \le 0.
 \ee  Notice that the integral on the right hand side is finite by a
 similar argument as for
 the case of $W_k(s)$.  So for fixed $s$, $W_k(s)=W(g, u, s \tau_k)$ is a
 non-increasing sequence of $k$.
 Using the lower bound on $W_k(s)$ in (\ref{wklowerb}), there exists a
 function
 $W_\infty(s)$ such that
 \[
 \lim_{k \to \infty} W_k(s) = \lim_{k \to \infty} W(g, u, s \tau_k) = W_\infty(s).
 \]

 Now we pick $s_0 \in [1, 2]$. Clearly we can find a subsequence $\{ \tau_{n_k} \}$,
  tending to infinity,
 such that
 \[
 W(g, u, s_0 \tau_{n_k})
 \ge W(g, u, (s_0+1) \tau_{n_k}) \ge W(g, u, s_0 \tau_{n_{k+1}}).
 \]Since
\[
\lim_{k \to \infty} W(g, u, s_0 \tau_{n_k})  =
\lim_{k \to \infty} W(g, u, s_0 \tau_{n_{k+1}})  = W_\infty(s_0),
\]we know that
\[
\lim_{k \to \infty} [ W(g, u, s_0 \tau_{n_k}) - W(g, u, (s_0+1) \tau_{n_k}) ] =0.
\]That is
\[
\lim_{k \to \infty} [ W_{n_k}(s_0) - W_{n_k}(s_0+1) ] =0.
\]Integrating (\ref{dwds}) from $s_0$ to $s_0 + 1$, we use the above to conclude that
\[
 \lim_{k \to \infty} \int^{s_0+1}_{s_0} \int s | Ric_{g_{n_k}} + Hess_{g_{n_k}} f_{n_k} - \frac{1}{2s}
 g_{n_k } |^2 u_{n_k} d\mu_{g_{n_k}(s)} ds = 0.
\]Therefore we have
\[
Ric_\infty + Hess_\infty f_\infty - \frac{1}{2s} g_\infty=0,
\]here
$f_\infty$ is defined by $(4\pi s)^{-n/2} e^{-f_\infty}
=u_\infty$. So the backward limit is a gradient shrinking Ricci
soliton.

Finally we need to show the soliton is non-flat.  Since the
original
 $\kappa$-solution is not a flat gradient shrinking soliton.
 Hence we know that
 \[
 W_k(s) < W_k(0) =W_0=0,
 \]where $W_0$ is the Euclidean $W$
entropy with respect to the standard Gaussian (see \cite{ni06}).
By (\ref{fk>c}), we know that the integrand in $W_k(s)$, $s \in
[1, 4]$, is bounded from below by a negative constant. Applying
Fatou's lemma on a sequence of exhausting domains, we find that
\be \lab{winftywk} \int \left[ s ( |\nabla f_\infty|^2 + R_\infty)
+f_\infty - n \right] u_\infty d\mu_{g_\infty(s)} \le W_k(s)
<W_0=0. \ee

If the gradient shrinking Ricci soliton $(\M, g_\infty)$ is flat,
it has to be $\R^n$. In fact, since
\[
 Hess_\infty f_\infty =\frac{1}{2s} g_\infty,
\]the universal cover  of $({\M}_\infty, g_\infty)$ is isometric to
$\R^n$.  Then it follows from a standard argument (for example,
see \cite[Page 203]{mt07}) that
$({\M}_\infty, g_\infty)$ is $\R^n$ with Euclidean metric.\\

 Note $\int_{M_\infty} u_{\infty} \le 1$ by Fatou's lemma again,
 we claim  that $$W(g_\infty, u_\infty, s) \doteq \int_{\R^n} \left[ s ( |\nabla f_\infty|^2 + R_\infty)
 +f_\infty - n
\right] u_\infty dx \ge 0.$$
  Let $\hat u=u_\infty/\Vert u_\infty
\Vert_1$. Recall that in the Euclidean space, the best constant of
the log Sobolev inequality is achieved by the Gaussian. Hence the
$W$-entropy associated with Gaussian is $0$ (see \cite{ni06}).
Therefore the $W$-entropy associated with $\hat u$ is nonnegative,
i.e.,
\[
W(g_\infty, \hat u, s) = \int_{\R^n} [ s \frac{|\nabla \hat
u|^2}{\hat u} -\hat u \ln \hat u - \frac{n}{2} (\ln 4 \pi s) \hat
u  - n  \hat u] dx \ge 0.
\]Now, $u_\infty = \Vert u_\infty \Vert_1  \hat u$ and
\[
W(g_\infty, u_\infty, s) \doteq \int_{\R^n} [ s \frac{|\nabla
u_\infty|^2}{ u_\infty} - u_\infty \ln  u_\infty - \frac{n}{2}
(\ln 4 \pi s) u_\infty - n  u_\infty] dx ,
\]this leads to that
\[
W(g_\infty, u_\infty, s) =\Vert u_\infty \Vert_1 W(g_\infty, \hat
u, s) -\Vert u_\infty \Vert_1 \ln \Vert u_\infty \Vert_1 \ge 0.
\] Hence the claim follows.\\

 The above claim is a contradiction of (\ref{winftywk}),
which means that $(\M_\infty, g_\infty)$ is not flat. This
finishes the proof of our theorem. \qed

\section{Appendix}
\label{app}

In this section we state (without proof) a uniform Sobolev inequality
under the Ricci flow, which is used in Sections 2 and 3. This result originally
appeared in \cite{z07} and \cite{z07err}, also see \cite{ye07} and \cite{hsu08}
for similar results in this direction.

\begin{theorem}
\label{thsobosmooth} Let $({\bf M}\sp n, g(t))$, $t \in [0, T_0)$, be a compact  solution to the Ricci flow (\ref{rf}), here $n \ge 3$ and $T_0
\le \infty$. Let $A_0$ and $B_0$ be two positive
numbers such that the following $L^2$ Sobolev inequality
holds for $({\bf M}, g(0))$: for all $v \in W^{1, 2}({\bf M})$,
\[
\left( \int_{\bf M} v^{2n/(n-2)} \,d\mu_{g(0)} \right)^{n/(n-2)} \le
A_0 \int_{\bf M} |\nabla v |^2 \,d\mu(g(0)) + B_0 \int_{\bf M} v^2
\,d\mu_{g(0)}.
\]Let $\lambda_0$ be the lowest eigenvalue of the operator $-4\Delta +R$ on $({\bf M}, g(0))$, i.e.,
\[
\lambda_0 = \inf_{ \Vert v \Vert_2 = 1} \int_{\M} ( 4 |\nabla v |^2
+ R v^2) d\mu_{g(0)}.
\]

Then the following conclusions are true.

(a). For all $t \in [0, T_0)$, there exist
positive functions $A(t)$ and  $B(t)$, depending only on the initial metric
$g(0)$ in terms of $A_0$, $B_0$, $\lambda_0$,  and $t$, such that, for all $v \in
W^{1, 2}({\bf M}, g(t))$, we have
\[
\bigg( \int v^{2n/(n-2)} d\mu_{g(t)} \bigg)^{(n-2)/n} \le A(t) \int (
|\nabla v |^2 + \frac{1}{4} R v^2 ) d\mu_{g(t)} + B(t)  \int v^2
d\mu_{g(t)},
\]here $R$ is the scalar curvature with respect to $g(t)$.
Moreover, if $R(x,0)>0, \forall x \in \M$, then $A(t)$ and $B(t)$
are constants independent of $t$.\\

(b). If the solution of Ricci flow is smooth for $t \in [0, 1)$
and becomes singular at $t=1$. Let $\tilde t = - \ln (1-t)$ and
$\tilde{g}(\tilde t) =\frac{1}{1-t} g(t)$, then under the
following Type I normalized Ricci flow
\[
\frac{\partial \tilde g}{\partial {\tilde t}} = - 2 \widetilde{Ric} + \tilde g,
\]
there exist positive constants $A$ and $B$, depending only on the
initial metric $g(0)$, such that, for all $v \in W^{1, 2}({\bf M},
\tilde g(\tilde t))$ and $\tilde t>0$, we have
\[
\bigg( \int v^{2n/(n-2)} d\mu_{\tilde g(\tilde t)} \bigg)^{(n-2)/n}
\le A \int ( |\tilde \nabla v |^2 + \frac{1}{4} \tilde R v^2 )
d\mu_{\tilde g(\tilde t)}+ B \int v^2 d\mu_{\tilde g(\tilde t)},
\]here $\tilde R$ is the scalar curvature with respect to $\tilde g(\tilde t)$.
\end{theorem}
\medskip

\begin{remark}  In fact, if for the initial metric $g(0)$ and $\forall v \in W^{1, 2}({\bf M}, g(0))$,
\[
\lambda_0 = \inf_{\Vert v \Vert_2 =1} \int(4 |\tilde \nabla v
|^2 + \tilde R v^2) d\mu_{g(0)}
>0,
\] then we also have $B=0$.
\end{remark}

The proof of the above theorem is based on the following uniform
log-Sobolev inequality: $\forall v \in W^{1, 2}({\bf M}, g(t))$
and $t>0$, then we have
\begin{equation}
\al
\lab{reslogsob1}
 \int_{\bf M}   v^2 \ln v^2 \,d\mu_{g(t)}
 \le \epsilon^2 \int_{\bf M} \big( 4 |\nabla  v |^2 +  R v^2 \big)
d\mu_{g(t)} -n \ln \epsilon + (t+\epsilon^2) \beta
+ \alpha,
\eal
\end{equation}
here $\alpha=\alpha(A_0, B_0, \lambda_0, n)$ and $\beta=\beta(A_0, B_0, \lambda_0, n)$ are positive constants.  Moreover $\beta=0$ provided $\lambda_0>0$.

For estimates of the constants $\alpha$ and $\beta$, please see \cite[Section 6.2]{bookCRC2010}.\\

{\bf Acknowledgement.}  This paper was essentially finished during Q.Z.'s visit at Cornell University. He wishes to thank Cornell Mathematics department for the kind hospitality.

%\bibliographystyle{plain}
%\bibliography{bio}
%\end{document}

%\begin{thebibliography}{10}
\def\cprime{$'$}

%\end{thebibliography}
\end{document}